\documentstyle[pb-diagram]{article}

\newcommand{\es}{\emptyset}

\newcommand{\ba}{\begin{array}}
\newcommand{\ea}{\end{array}}

\newtheorem{theorem}{Theorem}
\newtheorem{proposition}[theorem]{Proposition}
\newtheorem{lemma}[theorem]{Lemma}

\newtheorem{definition}[theorem]{Definition}
\newtheorem{corollary}[theorem]{Corollary}

\newcommand{\be}{\begin{enumerate}}
\newcommand{\ee}{\end{enumerate}}
\newcommand{\bi}{\begin{itemize}}
\newcommand{\ei}{\end{itemize}}
\newcommand{\bd}{\begin{description}}
\newcommand{\ed}{\end{description}}

\newcommand{\et}{\wedge}

\newcommand{\vel}{\vee}

\newcommand{\imp}{\rightarrow}

\newcommand{\beq}{\begin{eqnarray*}}
\newcommand{\eeq}{\end{eqnarray*}}
\newcommand{\seq}{\Rightarrow}

\author{ {F.Parlamento, F.Previale }
\\Department of Mathematics,  Computer Science and Physics
\\University of Udine,  via  Delle Scienze 206, 33100 Udine, Italy.
\\Department of Mathematics
\\University of Turin, via Carlo Alberto 10, 10123 Torino, Italy
\\e-mail: {\em franco.parlamento$@$uniud.it},  {\em flaviopreviale.cenasco@gmail.com}
}
 \title{Admissibility of the  Structural Rules in the Sequent Calculus with Equality 
\thanks{Work  partially supported  by the Italian PRIN  grant \emph{Mathematical Logic:  models, sets, computability}."
 }}
\date{}

\begin{document}


\maketitle




\noindent{{\bf Keywords :} Sequent Calculus, Structural Rules, Equality, Replacement Rules, Admissibility

\begin{abstract}

On the ground of the results in \cite{PP20} concerning the admissibility of the structural rules in sequent calculi with additional atomic rules,
 we develop a proof theoretic analysis 
for several extensions of the ${\bf G3[mic]}$ sequent calculi with rules for equality, including  the one  originally proposed by H.Wang  in the classic \cite{W60}.
In the classical case we relate our results with 
 the semantic tableau method for first order logic with equality. In particular 
we establish that, for languages without function symbols, in Fitting's alternative semantic tableau method in \cite{F96}
strictness (which does not allow  the repetition  of equalities  which are modified) can be imposed  together with the orientation of the replacement of equals.
 A significant 
progress is made toward extending that result to  languages with function symbols although whether that is possible or not remains to be settled.
We also briefly consider  systems that, in the classical case, are related to the 
semantic tableau method  in which one can expand branches by adding identities at will,
obtaining that also in that case strictness can be imposed.
Furthermore  we discuss to what extent  the strengthened form of the nonlengthening property of Orevkov obtained in  \cite{PP21} applies also to the present context.

\end{abstract}
\section{Introduction}
 
 In \cite{PP21} we have shown that  full cut elimination holds for the extension of  Gentzen's sequent calculi obtained by adding the Reflexivity Axiom $\seq t=t$, 
   and the 
 left introduction rules for $=$:
    
     \[
 \ba{clcl}
\Gamma\seq \Delta, F[x/r]& \vbox to 0pt{\hbox{$=_1$}}~~~&\Gamma\seq \Delta, F[x/r]&\vbox to 0pt{\hbox{$=_2$}}\\
 \overline{r=s,~ \Gamma\seq \Delta, F[x/s]}&&\overline{s=r, ~\Gamma\seq \Delta, F[x/s]}&
 \ea
 \]
where $F$ is a formula; $F[x/r]$ and $F[x/s]$, as in \cite{TS00}, denote the result of the replacement in $F$ of all free occurrences
of $x$ by $r$ or $s$ and $\Gamma$, $\Delta$ are finite multisets of formulae, with $|\Delta|=0$ in the intuitionistic case.
In \cite{PP21} the result is extended to other well motivated calculi with rules where $F[x/r]$ and $F[x/s]$ occur in the antecedent of the premiss and of  the conclusion. 
The purpose of this work is to introduce and study  corresponding sytems free of structural rules, some of which, in the classical case, are of particular interest in connection with
the semantic tableau method for first order logic with equality.
For that we have to refer to systems of that sort as far as logic is concerned such as    
the  multisuccedent systems for minimal, intuitionistic and classical logic originated with Dragalin's \cite{D88} and 
denoted by $\mbox{m-}{\bf G3[mic]}$ in  \cite{TS00}, that we will adopt as our logical systems.
 Since we will be dealing exclusively with such multisuccedent systems,  as remarked  in \cite{TS00} (pg. 83), the
 prefix $\mbox{m-}$  is redundant and we will drop it.
 Thus 
    ${\bf G3i}$ will 
 denote the multisuccedent  ${\bf G3}$ calculus for intuitionistic logic, ${\bf G3m}$ the analogous calculus for minimal logic, ${\bf G3c}$ the classical calculus 
  and  ${\bf G3[mic]}$ any of such three calculi.
We then adopt the Reflexivity Axiom in the form $\Gamma \seq t=t$, to be denoted by $\overline{\hbox{Ref}}$; 
restrict the formula $F$ in $=_1$ and $=_2$ to be atomic and, following the general pattern  exploited in  Kleene \cite{K52} to obtain sequent  calculi free of structural rules, we  repeat the principal formula $r=s$
in the antecedent of the premiss of the rules.
As we will show, that is both necessary and sufficient, and leads to what may be consider a most natural sequent calculus with equality  free of structural rules  and in the classical case coincides with the system first introduced, though semantic considerations, in the classic \cite{W60}.
We will denote with $\hbox{Rep}_1^r$ and $\hbox{Rep}_2^r$ the rules so obtained,  namely:

\[
\ba{clccl}
\underline{r=s, \Gamma \seq \Delta, P[x/r] }&\vbox to 0pt{\hbox{$\hbox{Rep}_1^{r}$}}&~~~~~~~~~&\underline{s=r, \Gamma\seq \Delta, P[x/r]}&\vbox to 0pt{$\hbox{Rep}_2^{r}$}\\
r=s,\Gamma\seq \Delta,  P[x/s] &&&s=r, \Gamma\seq \Delta, P[x/s]
\ea
\] 
where $P$ is atomic (possibly an equality), called the {\em context} formula,  while $r=s$ ($s=r$)  is called the {\em operating equality}  and $P[x/r]$ ($P[x/s]$) the {\em input} ({\em output}) formula.
In the classical case, the well known connection between such kind of calculi and the semantic tableau method for first logic with equality developed for example in \cite{J67} and \cite{F96},
add  motivations to those in \cite{PP21}, for the rules to follow:

\[
\ba{clccl}
\underline{r=s, P[x/r], \Gamma \seq \Delta }&\vbox to 0pt{\hbox{$\hbox{Rep}_1^{l}$}}&~~~~~~~~~&\underline{s=r,  P[x/r], \Gamma\seq \Delta}&\vbox to 0pt{$\hbox{Rep}_2^{l}$}\\
r=s,  P[x/s], \Gamma\seq \Delta &&&s=r, P[x/s], \Gamma\seq \Delta
\ea
\] which correspond to the tableau system in \cite{F96} pg 289 except that  {\em strictness} is required, namely
  the reuse of the formula in which the term replacement is operated in is not allowed, and, following (for Rep) the notation in \cite{TS00}), the rules:

\[
\ba{clccl}
\underline{r=s, P[x/s],P[x/r], \Gamma \seq \Delta}&\vbox to 0pt{\hbox{Rep'}}&~~~~&\underline{s=r,~ P[x/s],P[x/r], \Gamma \seq \Delta}&\vbox to 0pt{\hbox{Rep}}\\
r=s, P[x/s], \Gamma \seq \Delta&&~~&s=r,~ P[x/s], \Gamma \seq \Delta&
\ea
\] which  correspond to   the above  tableau  system in which strictness is not required.
In  such  tableau systems  a branch can be expanded by the addition of an identity $t=t$ at will. To that expansion rule it corresponds
the following Left Reflexivity Rule, denote by Ref in \cite{TS00}:
\[
\ba{cc}
\underline{t=t,\Gamma\seq \Delta}&\vbox to 0pt{\hbox{Ref}}\\
\Gamma\seq \Delta&
\ea
\]

 Our results will be based on the following fact that follows from the main result  in \cite{PP20}:
    for any set ${\cal R}$ of  atomic  rules  for equality that we will consider,
     if the structural rules are admissible in ${\cal R}$, identified with the calculus that consists of the initial sequents, including $\bot,\Gamma\seq \Delta$ in the intuitionistic and classical case,  and the rules  in ${\cal R}$,
     then they are admissible also in the calculus ${\bf G3[mic]}^{\cal R}$ obtained by adding the rules in ${\cal R}$ to ${\bf G3[mic]}$.

\section{Preliminaries on the logical calculi}

The sequent calculus  denoted by  ${\bf G3c}$  in \cite{TS00} (pg 83),  has the following initial sequents and rules, where $P$ is an atomic formula and $A, B$ stand for any formula in a first order language   (function symbols included) with bound variables distinct from the free ones,  and $\Gamma$ and $\Delta$ are finite multisets of formulae:

\

{\bf Initial sequents}
\[
P,\Gamma \seq \Delta, P
\]

{\bf Logical rules}

\[
\ba{clccl}
A,B, \Gamma \seq \Delta&\vbox to 0pt{\hbox{$L\et $}}&~~~~~~&\underline{\Gamma\seq \Delta, A~~~\Gamma\seq \Delta, B}& \vbox to 0pt{\hbox{$R\et$}}\\

\overline{A\et B, \Gamma \seq \Delta}&&&\Gamma \seq \Delta, A\et B\\

\\
\underline{A,\Gamma\seq \Delta~~~B,\Gamma \seq \Delta}& \vbox to 0pt{\hbox{$L\vel $}}&&\Gamma\seq \Delta, A, B& \vbox to 0pt{\hbox{$R\vel$}}\\

A\vel B, \Gamma \seq \Delta&&&\overline{\Gamma\seq\Delta, A\vel B}\\

\underline{\Gamma\seq \Delta, A~~~B,\Gamma \seq \Delta}&\vbox to 0pt{\hbox{$L\imp $}}&&A,\Gamma \seq \Delta, B&\vbox to 0pt{\hbox{$R\imp $}}\\

A\imp B, \Gamma \seq \Delta&&&\overline{\Gamma \seq \Delta, A\imp B}&\\

\\

&\vbox to 0pt{\hbox{$L\bot $}}&&&\\
\overline{\bot, \Gamma \seq \Delta}&&&&\\

\ea
\]

\[
\ba{clccl}
\underline{A[x/t], \forall x A,\Gamma\seq \Delta}&\vbox to 0pt{\hbox{$L\forall $}}&&\underline{\Gamma\seq \Delta, A[x/a]}&\vbox to 0pt{\hbox{$R\forall $}}\\

\forall xA, \Gamma\seq \Delta&&&\Gamma\seq \Delta, \forall x A\\

\\

\\
\underline{A[x/a], \Gamma\seq \Delta}&\vbox to 0pt{\hbox{$L\exists $}}&&\underline{\Gamma\seq \Delta, \exists x A,  A[x/t]}&\vbox to 0pt{\hbox{$R\exists $}}\\

\exists x A, \Gamma\seq \Delta&&&\Gamma\seq \Delta, \exists x A

\ea
\]

In ${\bf G3i}$ the rules $L\imp$, $R\imp$ and $R\forall$ are replaced by:

\[
\ba{clccl}
\underline{A\imp B, \Gamma \seq \Delta, A~~~B,\Gamma \seq \Delta}&\vbox to 0pt{\hbox{$L^i\imp $}}&~~~~~~~~&A,\Gamma \seq B&\vbox to 0pt{\hbox{$R^i\imp $}}\\

A\imp B, \Gamma\seq \Delta&&&\overline{\Gamma\seq\Delta, A\imp B}\\

\\
\\

&&&\Gamma\seq A[x/a]&\vbox to 0pt{\hbox{$R^i\forall $}}\\

&&&\overline{\Gamma\seq \Delta, \forall x A}&

\ea
\]

Finally ${\bf G3m}$ is obtained from ${\bf G3i}$ by replacing $L\bot$ by the initial sequents $\bot, \Gamma \seq \Delta, \bot$.

\

In all such systems  $a$ is a free variable that does not occur in the conclusion of $L\exists$ and $R\forall$.

\ 

${\bf G3[mic]}$ denotes any of the systems ${\bf G3m}$,   ${\bf G3i}$ or  ${\bf G3c}$.

 \ 

 The left and right weakening rules, LW and RW have the form:

\[
\ba{clccl}
\Gamma\seq \Delta&\vbox to 0pt{\hbox{LW}}&~~~~&\Gamma\seq \Delta&\vbox to 0pt{\hbox{RW}}\\
\cline{1-1}\cline{4-4}
A,\Gamma\seq \Delta&&&\Gamma\seq\Delta,A
\ea
\]

The left and right contraction rules, LC and RC have the form:

\[
\ba{clccl}
A.A, \Gamma\seq \Delta&\vbox to 0pt{\hbox{LC}}&~~~~&\Gamma\seq \Delta, A,A&\vbox to 0pt{\hbox{RC}}\\
\cline{1-1}\cline{4-4}
 A,\Gamma\seq \Delta&&&\Gamma\seq\Delta, A
 \ea
\]

$\hbox{LC}^=$ is the rule LC  in which the contracted formula $A$ is an equality.

The cut rule has  the form:
\[
\ba{clcclccl}
\Gamma\seq\Delta, A~~~~A,\Lambda \seq \Theta&\vbox to 0pt{\hbox{Cut }}&~~&&&~~&&\\
\cline{1-1}
\Gamma,\Lambda \seq   \Delta,\Theta&&&&&&\\
\ea
\]

Weakening, contraction and cut are the structural rules whose admissibility we are going to investigate.

In consequence of the more general result concerning the addition of atomic rules to the above sequent calculi established  in \cite{PP20}, for any set ${\cal R}$ of the above equality rules and  the further single premiss equality  rules to be introduced  in the sequel we have the following:  
\begin{theorem}\label{basic}[Theorem 1 in \cite{PP20}]
If the structural rules are admissible in ${\cal R}$, then they are admissible in  ${\bf G3[mic]}^{\cal R}$ as well.
\end{theorem}
that will be instrumental for the present work.

\ 

A further rule that will play an important  auxiliary role  is the following {\em congruence rule}:
\[
\ba{cl}
\underline{\Gamma_1\seq \Delta_1, r=s~~~~~~\Gamma_2\seq \Delta_2, P[x/r]}&\vbox to 0pt{\hbox{$\hbox{CNG}$}}\\
\Gamma_1,\Gamma_2\seq \Delta_1,\Delta_2, P[x/s]
\ea
\]

 {\bf Note}  The rule $CNG$ is among those used  in the extension of the system CERES in \cite{BL11}, pg.170.

\subsection{Admissibility of Weakening and Right Contraction}

The weakening rules are clearly   height preserving admissible in  the systems consisting of $\overline{\hbox{Ref}}$ and some of the equality rules.
The single premiss equality rules   modify at most one formula in the succedent of their  premiss.
Furthermore the initial sequents and those in  $\overline{\hbox{Ref}}$ remain initial sequents or in  $\overline{\hbox{Ref}}$ if all the formulae in their succedent, except the principal one,  are eliminated. By a straightforward induction on the height of derivations it follows that if $\Gamma\seq \Delta$ has a derivation  in the systems we are considering, then there  is a formula $A$ in $\Delta$ such that $\Gamma\seq A$ has a derivation of the same height.
That is the case also for the two premisses rule $\hbox{CNG}$ that eliminates a formula  from the succedent of its  first  premiss and modifies a single formula of the succedent of the second. As a consequence  the right contraction rule is  height-preserving admissible in all the  systems we are going to deal with.

\subsection{Basic equivalence theorem}

A basic tool for our investigation  is provided by the following proposition, where by an equality rule we mean any of the rules presented in the introduction
other than Ref and $\overline{\hbox{Ref}}$:

 \begin{proposition}\label{Equivalence}
  All the equality rules are equivalent in $\{\overline{{\em Ref}}$, {\em Cut, LC}$\}$ and $\{${\em Ref, Cut, LC}$\}$.
  \end{proposition}
  
  {\bf Proof}   $\overline{\mbox{ Ref}}$ is immediately derivable from Ref applied to the initial sequent $t=t\seq t=t$. Conversely
  Ref is derivable by applying the cut rule to its premiss $t=t,\Gamma\seq \Delta$  and the 
   the instance $\seq t=t$ of $\overline{\mbox{ Ref}}$. Therefore $\{ \overline{\hbox{Ref}},\hbox{Cut}, \hbox{LC}\}$ and  $\{ \hbox{Ref},\hbox{Cut}, \hbox{LC}\}$ are equivalent and it suffices to prove the  equivalence of the various rules with respect to one or the other of these two systems. We first show that if we  add  any one of the equality rules
  to such systems, then the following rule of Left Symmetry becomes derivable:
  \[
\ba{cc}
\underline{r=s,\Gamma\seq \Delta}&\vbox to 0pt{\hbox{Symm}}\\
s=r, \Gamma\seq \Delta&
\ea
\] 
Case 1.1. The rule added is  $\hbox{Rep}_1^r$. Then we have the following derivation of Symm:

   \[
   \ba{clcl}
s=r\seq s=s & ~~~&&\\
   \cline{1-1}
   s=r\seq r=s &&r=s, \Gamma\seq \Delta&\vbox to 0pt{$\hbox{Cut}$}\\
   \cline{1-3}
   \multicolumn{3}{c}{s=r,\Gamma\seq \Delta}
   
   \ea
   \]
Case 1.2. The rule added is $\hbox{Rep}_2^r$. Similar to Case 1.1

Case 2.1.  The rule added is  $\hbox{Rep}_1^l$. Then  we have the following derivation:

  \[
  \ba{cl}
  r=s, \Gamma\seq \Delta& \vbox to 0pt{$\hbox{LW}$}\\
  \cline{1-1}
  r=s, ~r=r, \Gamma\seq \Delta & \vbox to 0pt{$\hbox{Rep}_1^l$}\\
  \cline{1-1}
   r=s,~ s=r, \Gamma\seq \Delta & \vbox to 0pt{$\hbox{Rep}_1^l$}\\
   \cline{1-1}
    r=r, ~s= r, \Gamma\seq \Delta & \vbox to 0pt{$\hbox{Ref}$}\\
    \cline{1-1}
    s= r, \Gamma\seq \Delta & 
  \ea
  \] 
  Case 2.2. The rule added is  $\hbox{Rep}_2^l$. Then  the derivation is the same as for case 2.1, except that LW introduces $s=s$ and  $\hbox{Rep}_2^l$ is used
   instead of $\hbox{Rep}_1^l$.
  
  Case 3.1. The rule added is  Rep. Then  we have the following derivation:

   \[
   \ba{cl}
    \underline{ r=s, \Gamma \seq \Delta}&\vbox to 0pt{$\hbox{LW}$}\\
  \underline{ r=s, s=s, s=r, \Gamma \seq \Delta}&\vbox to 0pt{$\hbox{Rep}$}\\

\underline{ s=s, s=r, \Gamma \seq \Delta}&\vbox to 0pt{$\hbox{Ref}$}\\
  
   s=r, \Gamma \seq \Delta
   \ea
   \] 
   Case 3.2  The rule added is Rep'. Similar to Case 3..

   Case 4.   The rule added is CNG. Then we have the following derivation:
   
   \[
   \ba{clcl}
s=r\seq s=r~~~~\seq s=s& &&\\
   \cline{1-1}
   s=r\seq r=s &&r=s, \Gamma\seq \Delta&\vbox to 0pt{$\hbox{Cut}$}\\
   \cline{1-3}
   \multicolumn{3}{c}{s=r,\Gamma\seq \Delta}
   
   \ea
   \]
   Clearly the derivability of Symm makes equivalent the rules of the same type with index 1 and 2.
   Thus it suffices to  verify the equivalence (that does not depend on the availability of Symm)  between $\hbox{Rep}_1^r$ and  $\hbox{Rep}_2^l$; $\hbox{Rep}_1^l$ and  $\hbox{Rep}$; 
   $\hbox{Rep}_1^r$ and CNG.   We leave the easy details to the reader. $\Box$.

  \begin{corollary}\label{equivalentsystems}

All the systems ${\bf G3[mic]}^{\cal R}$, for  ${\cal R}$  that consists  of 
 $\overline{{\em Ref}}$ or ${\em Ref}$ and  of some of the  equality rules and such that the structural rules are admissible in ${\cal R}$,  are equivalent.

 \end{corollary}

 \section{Admissibility of the structural rules in systems  based on the Reflexivity Axiom}

 \subsection{Necessity of the repetition  of the operating equalities in the premiss of the equality rules }

We show that, as stated in the introduction,  the   addition of $\overline{\hbox{Ref}}$,  $~=_1$ and $=_2$  to ${\bf G3[mic]}$ is not sufficient to yield appropriate extensions free of structural rules. Actually 
 even if, beside $\overline{\hbox{Ref}}$,  $~=_1$ and $=_2$,  also the Cut rule is added, 
in the resulting system the left contraction rule  remains not admissible.

Let ${\cal R}=\{ \overline{Ref}, =_1, =_2, \hbox{Cut}\}$. We will prove that $LC$ is not admisssible in ${\bf G3c}^{\cal R}$ by showing that the following sequent:
\[
 *)~~a=f(a)\seq a=f(f(a))
\]
whose expansion $a=f(a), a=f(a) \seq a=f(f(a))$ is immediately derivable by means of an $=_2$-inference applied to $a=f(a)\seq a=f(a)$, is not derivable in ${\cal R}$.
In fact if $*)$ were derivable in ${\bf G3c}^{\cal R}$ (as for Proposition 7 in \cite{PP20}) $*)$ would have a derivation in which no logical inference different from $L\bot$ precedes a $=_1$, $=_2$ or Cut-inference. As a consequence $*)$  would be derivable in ${\cal R}$ itself,
which  is impossible.

In order to show that  $*)$ is not derivable in ${\cal R}$, we 
  first note that if a sequent $\Gamma\seq r=s$ is derivable in ${\cal R} $, then the sequent $\Gamma_=\seq r=s$, where $\Gamma_=$ denotes the multiset of equalities in $\Gamma$, has a derivation in ${\cal R} $ that involves only equalities. An easy induction on the height of such derivations shows that if
$\Gamma$ is a multiset of {\em identities} i..e equalities of the form $r=r$ and 
$\Gamma \seq r=s$ is derivable in ${\cal R}$ then $r=s$ is itself an identity ($r\equiv s$). That being noted, we prove the  following:

\begin{proposition} If $\Gamma$ is a multiset of identities and 
$E,\Gamma\seq E'$ is derivable in ${\cal R}$, where $E'$ coincides with 
$a=f(f(a))$ or with $f(f(a))=a$, then also $E$ has one of such two forms.
Hence
 $a=f(a) \seq a=f(f(a))$ is not derivable in ${\cal R}$.
\end{proposition}

{\bf Proof} We proceed by induction on the height of a derivation ${\cal D}$ in ${\cal R}$ of  $E,\Gamma\seq E'$.

If $h({\cal D})=0$, then $E,\Gamma\seq E'$ must be an initial sequent and $E$ coincides with $E'$ so that the claim is trivial.

 If $h({\cal D})>0$ and ${\cal D}$ ends with an $=_1$ inference that introduces $E$ in the antecedent, then ${\cal D}$ has the form:

\[
\ba{c}
{\cal D}_0\\
\Gamma \seq r=s\\
\cline{1-1}
E, \Gamma \seq E'
\ea
\]

By the previous remark $r=s$ is an identity $r=r$ and we note that the only possibilities of obtaining $E'$ by a substitution applied to $r=r$ is that $r$ coincides with  $a$ or with $f(f(a))$ in which  cases $E$ is necessarily $a=f(f(a))$ or $f(f(a))=a$.

 The same argument applies if ${\cal D}$ ends with an  $=_2 $-inference introducing $E$.
 
If  ${\cal D}$ ends with an  $=_1$ or $=_2 $-inference introducing  a formula in $\Gamma$, which is therefore an identity, the conclusion is a trivial consequence of the induction hypothesis.

If ${\cal D}$ ends with a Cut, we have two cases.

Case1. ${\cal D}$ has the form:
\[
\ba{ccc}
{\cal D}_0&~~~&{\cal D}_1\\
\Gamma_1 \seq A&&A,E,\Gamma_2\seq E'\\
\cline{1-3}

\multicolumn{3}{c}{E,\Gamma_1, \Gamma_2 \seq E'}

\ea
\]

In this case, looking at ${\cal  D}_0$  we have that $A$ is itself an identity so that it suffices to apply the induction hypothesis to ${\cal D}_1$ to conclude that $E$ is $a=f(f(a))$ or $f(f(a))= a $.

Case 2  ${\cal D}$ has the form:
\[
\ba{ccc}
{\cal D}_0&~~~&{\cal D}_1\\
E, \Gamma_1 \seq A&&A, \Gamma_2\seq E'\\
\cline{1-3}

\multicolumn{3}{c}{E,\Gamma_1, \Gamma_2 \seq E'}

\ea
\]
 By the induction hypothesis applied to ${\cal D}_1$  $A$ has one of the two forms
  $a=f(f(a))$
or $f(f(a))=a$ so that it suffices to apply the induction hypothesis to  ${\cal D}_0$ to conclude that the same holds for $E$.

That $a=f(a) \seq a=f(f(a))$ is not derivable in ${\cal R}$ follows by letting 
 $\Gamma$ be the empty set  and $E'$ the equality $a=f(f(a))$. 
  $\Box$

\section{Sufficiency of the repetition of the operating equalities in the premiss}

In this section we prove that the repetition of the operating equalities in the premiss of the $=_1$ and $=_2$-rules, which yields
 the $ \hbox{Rep}_1^r$ and  $ \hbox{Rep}_2^r$, suffices to yield a system, indeed a very natural one, for which the structural rules are admissible.

  \begin{theorem}\label{AdmissibilityR}
 For ${\cal R}_{12}^r =\{ \overline{{\em Ref}}, {\em Rep}_1^{r}, {\em Rep}_2^{r} \}$, the structural rules are admissible in ${\bf G3[mic]}^{{\cal R}_{12}^r}$.
  \end{theorem} 
 
{\bf Proof} By Theorem \ref{basic} it suffices to show that the structural rules are admissible in  ${\cal R}_{12}^r$. The admissibility of   $\hbox{LC}$   is straightforward, since the rules of ${\cal R}_{12}^r$ do not change the antecedent of their premiss. For the admissibility of  $\hbox{Cut}$ we transform a given derivation ${\cal D}$ in ${\cal R}_{12}^r + \hbox{Cut}$ into a derivation ${\cal D}'$
in $ \{ \overline{\hbox{Ref}}, \hbox{LC}, \hbox{CNG}, \hbox{Cut} \}$ by using  the following derivation of $\hbox{Rep}_1^r$  from CNG and $\hbox{LC}^=$:

\[
\ba{cl}
\underline{r=s\seq r=s~~~~r=s, \Gamma \seq \Delta, P[x/r]}&~~\vbox to 0pt{\hbox{CNG}}\\
\underline{ r=s, r=s, \Gamma \seq \Delta, P[x/s] }&\vbox to 0pt{{$\hbox{LC}^=$}}\\
r=s, \Gamma\seq \Delta, P[x/s]
\ea
\] and the derivation of $\hbox{Rep}_2^r$  from CNG and $\hbox{LC}^=$ that can be  obtained from that of $\hbox{Rep}_1^r$ 
thanks to the derivation of Symm from CNG shown
in  Case 4. of the proof of Proposition \ref{Equivalence}.

 From ${\cal D}'$  we  eliminate the applications of the $\hbox{Cut}$ rule in order  to obtain a derivation  ${\cal D}''$ in $ \{ \overline{\hbox{Ref}}, \hbox{LC}, \hbox{CNG} \}$.
 To show that this is possible, because of  the presence of the rule LC, we have to show that the following more general rule:
\[
\ba{c}
\Gamma\seq \Delta, A~~~~~A^n, \Lambda \seq \Theta\\
\cline{1-1}
\Gamma, \Lambda \seq \Delta, \Theta
\ea
\]
 where $A^n$ denotes the multiset that contains $A$  $n$ times and nothing else, is admissible in  $ \{ \overline{\hbox{Ref}}, \hbox{LC}, \hbox{CNG}\}$. That is shown by a straightforward induction on the height of the derivation of $A^n, \Lambda \seq \Theta$.

Then to obtain, from ${\cal D}''$,  the desired cut-free derivation in ${\cal R}_{12}^r$ of the endsequent of ${\cal D}$, it suffices to  exploit the admissibility of $\hbox{LC}$  and $\hbox{CNG}$ in ${\cal R}_{12}^r$.  
 The admissibility of  $\hbox{CNG}$ in ${\cal R}_{12}^r+ \hbox{LC}$, hence in   ${\cal R}_{12}^r$,   can be proved by induction on the height of the derivation  of its  first  premiss (see \cite{PP21} for the analogous result for the sequent calculi with structural rules).
In fact let  ${\cal D}$ be of the form:
 
 \[
 \ba{cccl}
 {\cal D}_0&&{\cal D}_1\\
 \Gamma' \seq \Delta, r=s&~~~&\Lambda\seq \Theta, P[x/r] &\vbox to 0pt{\hbox{CNG}} \\
 \cline{1-3}
 \multicolumn{3}{c}{\Gamma,\Lambda \seq \Delta, \Theta, P[x/s]}

 \ea
 \]where  ${\cal D}_0$ and ${\cal D}_1$ are derivations in ${\cal R}_{12}^r + \hbox{LC}$. 
 We have to show that also the conclusion of ${\cal D}$ is derivable in ${\cal R}_{12}^r + \hbox{LC}$.
 If $r$ and $s$ coincide, then the conclusion is obtained by weakening the conclusion of ${\cal D}_1$. Assuming $r$ is distinct from  $s$, we proceed by induction on the height $h({\cal D}_0)$ of ${\cal D}_0$.
 
  If $h({\cal D}_0)=0$ and 
${\cal D}_0$ is  an initial sequent with principal formula common to $\Gamma$ and $\Delta$, then the conclusion of ${\cal D}$ is also an initial sequent and the given of CNG-inference  can be eliminated,
while if it is of the form 
$r=s,\Gamma' \seq \Delta, r=s$, then ${\cal D}$, namely
 \[
 \ba{ccc}
 &&{\cal D}_1\\
 r=s,\Gamma'\seq \Delta, r=s&& \Lambda\seq \Theta, P[x/r]\\
 \cline{1-3}
 \multicolumn{3}{c}{ r=s,\Gamma',  \Lambda \seq \Delta, \Theta, P[x/s]}
 \ea
 \]
 is transformed into:

 \[
 \ba{cl}
 {\cal D}_1&\\
 \Lambda\seq \Theta, P[x/r]&\vbox to 0pt{$\hbox{LW}$}\\
  \overline{r=s,~\Gamma', \Lambda\seq \Delta, \Theta, P[x/r]}&\vbox to 0pt{$\hbox{Rep}_1^r$}\\\
  \overline{ r=s,\Gamma',  \Lambda \seq \Delta, \Theta, P[x/s]}&
 \ea
 \]
 If $h({\cal D}_0)>0$ and ${\cal D}_0$ ends with an $\hbox{Rep}_1^r$- inference and the principal formula occurs in $\Delta$ then the derivation of the conclusion is obtained as a straightforward consequence of the induction hypothesis. On the other hand if the principal formula is $r=s$ of the form $r^\circ[x/q]= s^\circ[x/q]$, with  ${\cal D}$ 
of  the form:
 \[
 \ba{ccc}
 {\cal D}_{00}&&\\
\underline{ p=q, \Gamma'\seq \Delta, r^\circ[x/p]=s^\circ[x/p]}&~~~~&{\cal D}_1\\
 p=q, \Gamma'\seq \Delta,  , r^\circ[x/q]=s^\circ[x/q]&&\Lambda\seq \Theta, P[x/r^\circ[x/q]]\\
 \cline{1-3}
 \multicolumn{3}{c}{ p=q,\Gamma',  \Lambda \seq \Delta, \Theta, P[x/s^\circ[x/q]]}
  \ea
 \]
${\cal D}$ can be      transformed into:
 
 \[
 \ba{cccl}
 &&{\cal D}_1&\\
 &~~~& \Lambda\seq \Theta, P[x/r^\circ[x/q]]&\vbox to 0pt{$\hbox{LW}$}\\
 {\cal D}_{00}&~~~& \overline{p=q, \Lambda\seq \Theta, P[x/r^\circ[x/q]]}&\vbox to 0pt{$\hbox{Rep}_2^r$}\\
 p=q, \Gamma' \seq \Delta, r^\circ[x/p]=s^\circ[x/p]&~~~& \overline{p=q, \Lambda\seq \Theta, P[x/r^\circ[x/p]]}&\vbox to 0pt{$\hbox{ind}$}\\
 \cline{1-3}
  \multicolumn{3}{c}{ \underline{p=q,~p=q,~\Gamma',  \Lambda \seq \Delta, \Theta, P[x/s^\circ[x/p]]}}&\vbox to 0pt{$\hbox{Rep}_1^r$}\\

  \multicolumn{3}{c}{ \underline{p=q,~p=q,~\Gamma',  \Lambda \seq \Delta, \Theta, P[x/s^\circ[x/q]]}}&\vbox to 0pt{{$\hbox{LC}^=$}}\\  

  \multicolumn{3}{c}{ p=q,~\Gamma',  \Lambda \seq \Delta, \Theta, P[x/s^\circ[x/q]]}
  \ea
 \]where  $\hbox{ind}$ means that,  by induction hypothesis,  the given derivations in ${\cal R}_{12}^r+ \hbox{LC}$ of the sequents above the line can be transformed
  into a derivation in ${\cal R}_{12}^r+ \hbox{LC}$ of the sequent below the line.
 If the premiss is obtained by an $\hbox{Rep}_2^r$  the argument is the same except that in the transformed derivation  we use $\hbox{Rep}_1^r$ in place of $\hbox{Rep}_2^r$
 and conversely. The case in which the first premiss is obtained by means of an $\hbox{LC}$-inference
 is straightforward. $\Box$

 \begin{corollary} \label{sistema completo}

 ${\cal R}_{12}^{rl} =\{ \overline{{\em Ref}}, {\em Rep}_1^{l},  {\em Rep}_2^{l},  {\em Rep}_1^{r}, {\em Rep}_2^{r} \}$ and 

  ${\cal R}_{12}^r $  are equivalent systems over which the structural rules are admissible.
 
 \end{corollary}
 
 {\bf Proof} Obviously ${\cal R}_{12}^r $ is a subsystem of  ${\cal R}_{12}^{rl}$. The converse holds by the previous Theorem and 
 the equivalence of  the equality rules over systems containing $\{\overline{\hbox{ Ref}}, \hbox{ Cut}, \hbox{ LC}\}$ established in
 Proposition \ref{Equivalence}. $\Box$

 \

Theorem \ref{AdmissibilityR} can be strengthened by requiring that, when the context formula is an equality,  the rules $\hbox{Rep}_1^r$ and $\hbox{Rep}_2^r$
change only its right-hand side.  Let  $\hbox{Rep}_1^{r =_r}$ and $\hbox{Rep}_2^{r=_r}$ be the restrictions of $\hbox{Rep}_1^r$ and $\hbox{Rep}_2^r$ obtained in that way.

   \begin{theorem}\label{AdmissibilityRr}
 The system  ${\cal R}_{12}^{r=_r} =\{ \overline{{\em Ref}}, {\em Rep}_1^{r=_r}, {\em Rep}_2^{r=_r} \}$ is equivalent  to  ${\cal R}_{12}^r$, hence 
  the structural rules are admissible in ${\bf G3[mic]}^{{\cal R}_{12}^{r=_r}}$.
  \end{theorem}

 {\bf Proof} It suffices to show that if a sequent of the form $\Gamma \seq p=q$ is derivable in ${\cal R}_{12}^r$, then it is derivable in ${\cal R}_{12}^{r=_r}$ as well.
 Given  a derivation ${\cal D}$ in ${\cal R}_{12}^r$ of $\Gamma\seq p=q$ we proceed by induction on the number of $\hbox{Rep}^r_1$ or $\hbox{Rep}^r_2$ -inferences that act on an equality but are not  $\hbox{Rep}^{r=_r}_1$ or $\hbox{Rep}^{r=_r}_2$ -inferences, to be called  {\em undesired}  inferences.
  If there are none we are done. Otherwise we select the topmost one call it $J$. Let us assume that it is of the form:
 
 \[
 \ba{cl}
r=s,  \Gamma^- \seq p'^{\circ}[x/r]= q'& \vbox to 0pt{$\hbox{Rep}_1^r$}\\
 \cline{1-1}
 r=s,  \Gamma^- \seq p'^{\circ}[x/s]= q'
 \ea
 \]
Since an initial sequent of the form $t=t', ~\Gamma \seq t = t'$ is derivable from  $t=t', \Gamma \seq t=t$ by means of a $\hbox{Rep}_1^{r=_r}$-inference,
  we may assume that the initial sequent of ${\cal D}$ has the form
 \[
 r=s, \Gamma^-\seq p'^{\circ}[x/r] = p'^{\circ}[x/r]
 \]
 If we  replace the initial sequent of ${\cal D}$ by:
 
\[
\ba{cl}
 r=s, \Gamma^- \seq p'^{\circ}[x/s] =  p'^\circ[x/s] & \vbox to 0pt{$\hbox{Rep}_2^{r=_r}$}\\
 \cline{1-1}
  r=s, \Gamma^- \seq p'^{\circ}[x/s] =  p'^\circ[x/r] 
\ea
\]and the successive left-hand sides $ p'^\circ[x/r]$ of the right equalities of ${\cal D}$ down to the premiss of $J$ by $ p'^\circ[x/s]$ we obtain the conclusion of $J$
that therefore can be eliminated from the given derivation of $\Gamma \seq p=q$ thus obtaining a derivation that has one less undesired inference than ${\cal D}$.
If  $J$ is  an $\hbox{Rep}^r_2$ the argument is the same except that the new 
initial  inference is a $\hbox{Rep}^{r=_r}_1$-inference rather than a  $\hbox{Rep}^{r=_r}_2$-inference.  $\Box$

 \section{Limiting the scope of  replacement}
  
 Let   $\hbox{Rep}_1^{r/=}$ and $ \hbox{Rep}_2^{r/=} $ be the rules  $\hbox{Rep}_1^{r=_r}$ and $ \hbox{Rep}_2^{r=_r} $ restricted to context  formulae that are equalities 
  and $\hbox{Rep}_1^{l/(=)}$ and $ \hbox{Rep}_2^{l/(=)} $ be the rules  $\hbox{Rep}_1^{l}$ and $ \hbox{Rep}_2^{l} $ restricted to context  formulae that are not  equalities. 
 \

 \begin{theorem}\label{AdmissibilityAT}
  Let  ${\cal R}^{r/=}_{l/(=)}$ be $\{ \overline{{\em Ref}},   {\em Rep}_1^{~l/(=)}, 
 {\em Rep}_2^{~l/(=) }, ~{\em Rep}_1^{r/=}, {\em Rep}_2^{~r/=} \}$. 
 
  ${\cal R}^{r/=}_{l/(=)}$ is equivalent to ${\cal R}_{12}^r$,
  therefore    
   the structural rules are admissible in
    ${\bf G3[mic]}^{{\cal R}^{r/=}_{l/(=)}}$.   
 \end{theorem} 
 
 {\bf Proof} 
 By Corollary \ref{sistema completo} every sequent derivable in  ${\cal R}^{r/=}_{l/(=)}$  is derivable in ${\cal R}_{12}^r$ as well.
  For the converse we note that  if $\Gamma\seq \Delta$ is derivable in ${\cal R}_{12}^{r}$, then there is a formula $A$ in $\Delta$ such that $\Gamma\seq A$ is  also derivable in that system. If $A$ is an equality, then
 the derivation of $\Gamma \seq A$ can use only $\hbox{Rep}_1^{r/=}$ and $\hbox{Rep}_2^{r/=}$, so that it belongs to ${\cal R}^{r/=}_{l/(=)}$.
 If $A$ is not an equality  we proceed by induction on the height of  the   derivation ${\cal D}$ in   ${\cal R}_{12}^{r}$ of $\Gamma\seq A$ to show that it can be transformed into a derivation (of the same height) in ${\cal R}^{r/=}_{l/(=)}$. If $h({\cal D})=0$, then $\Gamma\seq A$ is  an initial sequent  and the conclusion is obvious. If $h({\cal D})= n+1$, then ${\cal D}$ ends either with an $ \hbox{Rep}_1^{r}$-inference  or with an $ \hbox{Rep}_2^{r}$-inference. Let us assume, for example, that 
 ${\cal D}$ ends with a $ \hbox{Rep}_1^{r}$-inference. Then ${\cal D}$  has the form:
 
 \[
 \ba{c}
 P_1, \Gamma_1 \seq P_1\\
 {\cal D}_0\\
 \underline{r=s, \Gamma_n \seq P_n[x/r]}\\
 r=s, \Gamma_n\seq P_n[x/s]
 \ea
 \]
  By induction hypothesis there is a derivation ${\cal D}_0'$ in ${\cal R}^{r/=}_{l/(=)}$  (of height $n$) of 
  
  $r=s, \Gamma_n \seq P_n[x/r]$. ${\cal D}_0'$ has the form:
  
  \[
  \ba{c}
  \underline{r=s, P_n[x/r], \Lambda\seq P_n[x/r]}\\
  r=s, \Lambda' \seq P_n[x/r]\\
  \vdots\\

  r=s, \Gamma_n \seq P_n[x/r]
   \ea
  \]
In fact  $\hbox{Rep}_1^{l/(=)}$ and $\hbox{Rep}_2^{l/(=)}$ do not introduce any new equality in their conclusion, so that
 all the equalities in the endsequent of ${\cal D}_0'$, in particular $r=s$,   are present in the
 antecedent of every sequent in ${\cal D}_0'$.  If we replace all the occurrences of $P_n[x/r]$ in the succedents of the sequents of ${\cal D}_0$
by $P[x/s]$  and introduce an initial ${\hbox{Rep}}_2^{l/(=)}$-inference replacing $s$ by $r$ in $P_n[x/r]$ we obtain the desired derivation ${\cal D}'$ 
in ${\cal R}^{r/=}_{l/(=)}$ (of height $n+1$), namely:
 
 \[
  \ba{c}
  \underline{r=s, P_n[x/s], \Lambda \seq P_n[x/s]}\\
  \underline{r=s, P_n[x/r], \Lambda\seq P_n[x/s]}\\
  r=s, \Lambda' \seq P_n[x/s]\\
  \vdots\\

  r=s, \Gamma_n \seq P_n[x/s]
   \ea
  \]$\Box$
 
 \

 Clearly the proof goes through without any  change if $\hbox{ Rep}_1^{r/=}$ and $ \hbox{ Rep}_2^{r/=} $ are restricted to   $\hbox{Rep}_1^{r/=_r}$ 
 and $ \hbox{Rep}_2^{r/=_r} $ that change only the right-hand side of the equality that they transform. 
 
 Thus, letting  ${\cal R}^{r/=_r}_{l/(=)}=
 \{ \overline{{\em Ref}},   {\em Rep}_1^{~l/(=)}, 
 {\em Rep}_2^{~l/(=) }, ~{\em Rep}_1^{r/=_r}, {\em Rep}_2^{~r/=_r} \}$,
we have the following stregthened form of the previous Theorem:

 \begin{theorem}\label{AdmissibilityATr}

  ${\cal R}^{r/=_r}_{l/(=)}$ is equivalent to ${\cal R}_{12}^r$,
  therefore    
   the structural rules are admissible in
    ${\bf G3[mic]}^{{\cal R}^{r/=_r}_{l/(=)}}$.   
 \end{theorem}


Interpreted in terms of the alternate tableau system in \cite{F96}, pg. 294  where  a   branch can be closed if the  negation of an identity $\neg t=t$ appears on it, and left-right and {\bf right-left}  replacement can be applied to atomic formulae
and to negation of equalities,
this result, in the classical case, means that, 
 strictness can be imposed (no reuse of formulae in which a replacement is performed is  allowed)
  and the replacement rule can be applied only to atomic formula that are not equalities and to the right-hand side of negation of equalities.

\section{Orienting replacement  in languages without function symbols}

 We   prove  that for languages free of function symbols the structural rules are admissible in ${\cal R}^{rl}_{2}$ by showing 
  that for such languages ${\cal R}^{rl}_{2}$ is in fact equivalent to  ${\cal R}^r_{12}$.
   The same holds, with the same proofs,   for ${\cal R}^{rl}_{1}$. 
   
   \ 
  
 {\bf  Notation} In the following $a, b, c$   will stand for constants or  free variables and 
  $a\approx b$ may denote  either one of $a=b$ or $b=a$.
 
 \begin{definition}
 A chain of equalities connecting $a$ and $b$ denoted by $\gamma(a,b)$ is a set of equalities that can be arranged into a sequence of the form
 $a\approx a_1, a_1\approx a_2, \ldots, a_{n-1}\approx b$.
 The empty set is a chain  that connects any term  with itself.
   \end{definition}
   
   \begin{lemma} \label{Chain derivability} Given a chain $\gamma(a,b)$ and an atomic formula $A$ with at most one occurrence of $x$
   
   \bi
   \item[a)]   $\gamma(a,b)\seq a=b$ is derivable in ${\cal R}^{rl}_2$
   \item[b)]  $A[x/a], \gamma(a,b) \seq A[x/b]$ is derivable in  ${\cal R}^{rl}_2$
   \ei
   
   \end{lemma}
  
  {\bf Proof} In both cases  we proceed by induction on the length $n$ of  $a\approx a_1, a_1\approx a_2, \ldots, a_{n-2}\approx a_{n-1}, a_{n-1}\approx  b$.
  
 $a)$  For $n=0$, $\gamma(a,b)=\es$ and  $a\equiv b$ so that
  $\gamma(a,b)\seq a=b$ is the  instance  $\seq a=a$ of $\overline{\hbox{Ref}}$. For $n=1$,   $\gamma(a, b)$ is either $a=b$ or $b=a$.
  In the former case $\gamma(a, b)\seq a=b$ is the  initial sequent  $a=b\seq a=b$, while in the latter case it has the following derivation in ${\cal R}^{rl}_2$:
  
  \[
  \ba{cl}
  \underline{b=a\seq a=a}& \vbox to 0pt{$\hbox{Rep}_2^r$}\\
  b=a\seq a=b
  \ea
  \]
  
  Assume $n>1$. If $a_{n-1} \approx b$ is $a_{n-1} = b$, by induction hypothesis:
  
  \[
  a\approx a_1, \ldots, a_{n-2}\approx b \seq  a=b
  \] has a derivation   in    ${\cal R}^{rl}_2$ from which we obtain the desired  derivation  in    ${\cal R}^{rl}_2$ by the admissibility of LW
  that allows for the introduction of $a_{n-1} = b$ and a $\hbox{Rep}_2^l$-inference using $a_{n-1} = b$ as operating equality, namely:
  
  \[
  \ba{ll}

   \underline{a\approx a_1, \ldots, a_{n-2}\approx b \seq  a=b }& \vbox to 0pt{\hbox{LW}}\\
   \underline{ a\approx a_1, \ldots, a_{n-2}\approx b,~ a_{n-1}= b \seq  a=b} & \vbox to 0pt{\hbox{$\hbox{Rep}_2^l$}}\\
    a\approx a_1, \ldots, a_{n-2}\approx a_{n-1},~ a_{n-1}= b \seq  a=b
  \ea
  \]
   If $a_{n-1} \approx b$ is $b =a_{n-1} $, by induction hypothesis:
   
   \[
     a\approx a_1, \ldots, a_{n-2}\approx a_{n-1} \seq  a=a_{n-1}
   \]
   has a derivation ${\cal D}$  in    ${\cal R}^{rl}_2$ from which we obtain the desired  derivation  in    ${\cal R}^{rl}_2$ by the admissibility of LW
  that allows for the introduction of $b= a_{n-1} $ and a $\hbox{Rep}_2^r$-inference using $b= a_{n-1}$, namely:

  \[
  \ba{ll}

   \underline{a\approx a_1, \ldots, a_{n-2}\approx a_{n- 1} \seq  a=a_{n-1} }& \vbox to 0pt{\hbox{LW}}\\
   \underline{ a\approx a_1, \ldots, a_{n-2}\approx a_{n-1},~  b= a_{n-1} \seq  a=a_{n-1}} & \vbox to  0pt{\hbox{$\hbox{Rep}_2^r$}}\\
    a\approx a_1, \ldots, a_{n-2}\approx a_{n-1}, ~b= a_{n-1} \seq  a=b
  \ea
  \]
  
  $b)$ For $n=0$, $A[x/a], \gamma(a,b)\seq A[x/b])$ reduces to the initial sequent 
  
  $A[x/a]\seq A[x/a]$. For $n=1$ we have the following derivations,
  depending on whether $a\approx b$ is $a=b$ or $b=a$:
  
   \[
  \ba{clcl}
   \underline{A[x/b],~a=b \seq A[x/b]}&  \vbox to 0pt{\hbox{$\hbox{Rep}_2^l$}}~~~~~~&\underline{A[x/a], b=a \seq A[x/a]}&\vbox to 0pt{$\hbox{Rep}_2^r$}\\
 A[x/a], ~ a=b \seq A[x/b]&&A[x/a] ~b=a \seq A[x/b]
  \ea
  \]

For $n>1$ the argument is similar to that in $a)$. If  $a_{n-1}\approx b$ is $a_{n-1}=b$, we note that by induction hypothesis we have a derivation in ${\cal R}^{rl}_2$ 
of 
\[
 A[x/a],   a\approx a_1, \ldots, a_{n-2}\approx b \seq   A[x/b]
\]
  from which the desired derivation is obtained by a weakening introducing   $a_{n-1}=b$ followed by  a $\hbox{Rep}_2^l$-inference transforming  $a_{n-2}\approx b$
  into   $a_{n-2}\approx a_{n-1}$.
  
  If  $a_{n-1}\approx b$ is $b= a_{n-1}$, by induction hypothesis we have a derivation in ${\cal R}^{rl}_2$ 
of 
\[
 A[x/a],   a\approx a_1, \ldots, a_{n-2}\approx a_{n-1} \seq A[x/a_{n-1}]
\]
  from which the desired derivation is obtained by a weakening introducing  $b= a_{n-1}$ and a $\hbox{Rep}_2^r$-inference transforming  $A[x/a_{n-1}]$
  into   $ A[x/b] $. $\Box$

  \begin{lemma} \label{derivazione catena}
  Given an atomic formula $A$,  $m$ variables $x_1,\ldots, x_m$ having at most one occurrence in $A$
  and $m$ chains $\gamma_1(a_1, b_1), \ldots, \gamma_m(a_m, b_m)$ the sequent:
  \[
  A[x_1/a_1, \ldots, x_m/ a_m], \gamma_1(a_1, b_1), \ldots, \gamma_m(a_m, b_m) \seq A[x_1/b_1, \ldots, x_m/ b_m]
  \]is derivable in ${\cal R}^{rl}_2$.
   \end{lemma}
   
   {\bf Proof} We proceed by a principal induction  on $m$ and a secondary induction on the length of $\gamma_m(a_m, b_m)$.
   For $m=1$ the claim reduces to the previous lemma part $b)$. Assuming it holds for $m-1$ we have
    \[
    \ba{r}
1)~~ A[x_1/a_1, \ldots, x_{m-1}/a_{m-1}, x_m / a_m), \gamma_1(a_1, b_1), \ldots, \gamma_{m-1} (a_{m-1}, b_{m-1})\seq \\
 \seq A[x_1/b_1, \ldots, x_{m-1}/ b_{m-1}, x_m/a_m]
 \ea
  \] as well as  
     \[
    \ba{r}
2)~~ A[x_1/a_1, \ldots, x_{m-1}/a_{m-1}, x_m / b_m], \gamma_1(a_1, b_1), \ldots, \gamma_{m-1} (a_{m-1}, b_{m-1})\seq \\
 \seq A[x_1/b_1, \ldots, x_{m-1}/ b_{m-1}, x_m/ b_m]
 \ea
  \] 
  are derivable in ${\cal R}^{rl}_2$.
  Then we can proceed  by induction on the length $l$ of $\gamma_m(a_m, b_m)$ to show that also
   \[
  A[x_1/a_1, \ldots, x_m/ a_m], \gamma_1(a_1, b_1), \ldots, \gamma_m(a_m, b_m) \seq A[x_1/b_1, \ldots, x_m/ b_m]
  \]is derivable in ${\cal R}^{rl}_2$.
  If $l=0$ then $a_m\equiv b_m$ and the conclusion is immediate. If $l =1$ then 
  $\gamma_m(a_m, b_m)$ is either $a_m=b_m$ or $b_m=a_m$. In the first case we weaken the sequent $2)$ by adding $a_m=b_m$ and then apply a 
   $\hbox{Rep}_2^l$-inference to tranform $ A[x_1/a_1, \ldots, x_{m-1}/a_{m-1}, x_m / b_m]$ in the antecedent of $2)$ into $ A[x_1/a_1, \ldots, x_{m-1}, x_m / a_m]$.
   Similarly if   $\gamma_m(a_m, b_m)$  is $b_m=a_m$, we add  $b_m=a_m$ to the antecedent of $1)$ and then apply a $\hbox{Rep}_2^r$-inference 
   to transform   $A[x_1/b_1, \ldots, x_{m-1}/ b_{m-1}, x_m / a_m]$ in the consequent of $2)$ into  $A[x_1/b_1, \ldots, x_{m-1}/ b_{m-1}, x_m/ b_m]$.
  For  $l >1$ let $\gamma(a_m, b_m)$ be $a_m\approx a_m^1, a_m^1\approx a_m^2, \ldots, a_m^{l-2}\approx a_m^{l-1},   a_m^{l-1}\approx b_m$.
  If $a_m^{l-1}\approx b_m$ is $b_m= a_m^{l-1}$ we note that by induction hypothesis:

    \[
    \ba{l}
 A[x_1/a_1, \ldots, x_{m-1}/a_{m-1}, x_m / a_m], \gamma_1(a_1, b_1), \ldots, \gamma_{m-1}(a_{m-1}, b_{m-1}),\\
  a_m\approx a_m^1, a_m^1\approx a_m^2, \ldots, a_m^{l-2}\approx a_m^{l-1}  \seq  A[x_1/b_1, \ldots, x_{m-1}/ b_{m-1}, x_m/a_m^{l-1}]
 \ea
  \] 
is derivable in  ${\cal R}^{rl}_2$. Then it suffices to weaken the antecedent by adding $b_m= a_m^{l-1}$ and apply a  $\hbox{Rep}_2^r$-inference to transform 
$ A[x_1/b_1, \ldots, x_{m-1}/ b_{m-1}, x_m/a_m^{l-1}]$ into  $ A[x_1/b_1, \ldots, x_{m-1}/ b_{m-1},x_m/ b_m]$ to obtain the desired derivation  in  ${\cal R}^{rl}_2$ of 

  \[
  \ba{l}
*)~~~ A[x_1/a_1, \ldots, x_{m-1}/a_{m-1}, x_m / a_m], \gamma_1(a_1, b_1), \ldots, \gamma_{m-1}(a_{m-1}, b_{m-1}),\\
 a_m\approx a_m^1, a_m^1\approx a_m^2, \ldots, a_m^{l-2}\approx a_m^{l-1},     b_m= a_m^{l-1} \seq A[x_1/b_1, \ldots, x_{m-1}/ b_{m-1}, x_m/b_m]
\ea
  \] 
On the other hand if  $a_m^{l-1}\approx b_m$ is $a_m^{l-1}=b_m$ we note that by induction hypothesis  there is  a derivation in  ${\cal R}^{rl}_2$ of

   \[
   \ba{r}
 A[x_1/a_1, \ldots, x_{m-1}/a_{m-1}, x_m / a_m], \gamma_1(a_1, b_1), \ldots, \gamma_{m-1}, a_m\approx a_m^1, a_m^1\approx a_m^2, \ldots, a_m^{l-2}\approx b_m  \seq\\   \seq A[x_1/b_1, \ldots, x_{m-1}/ b_{m-1}, x_m/b_m]
 \ea
  \] 
  that can be weakened by the addition of $a_m^{l-1}=b_m$ in the antecedent to be used to transform, by means of a  $\hbox{Rep}_2^l$-inference,  $ a_m^{l-2}\approx b_m $ into   $ a_m^{l-2}\approx a_m^{l-1} $ in order to obtain 
  a derivation of $*)$ in  ${\cal R}^{rl}_2$.
  $\Box$

  \begin{lemma} \label{chain in}
  
  \bi
  \item[a)] If $\Gamma \seq a=b$ is derivable in ${\cal R}_{12}^r$, then $\Gamma$ includes a chain $\gamma(a,b)$.
  
  \item[b)] If $A$ is not an equality and $\Gamma \seq A$ is derivable in ${\cal R}_{12}^r$, then for some $m$ 
  there are  two $m$-tuples
   $a_1,\ldots  a_m$ and $b_1, \ldots, b_m$,  such that $A$ has the form $A^\circ[x_1/b_1, \ldots, x_m/b_m]$  and 
$\Gamma $ contains $A^\circ[x_1/a_1, \ldots, x_m/a_m]$ as well as  $m$ chains $\gamma_1(a_1, b_1), \ldots, \gamma_m(a_m, b_m)$. 
  \ei

  \end{lemma}\label{catena in}
  
  {\bf Proof} By Theorem 
  \ref{AdmissibilityRr} 
  we can  proceed by induction on the height of a derivation ${\cal D}$ in ${\cal R}_{12}^{r=_r}$ of $\Gamma \seq a=b$ or $\Gamma \seq A$. 
  
  a) If $h({\cal D})=0$ then  $\Gamma\seq a=b$ is an instance of $\overline{\hbox{Ref}}$ i.e. $a\equiv b$ and we can let  $\gamma(a,b)=\es$ or
   it is an initial sequent, i.e.
  $a=b $ occurs in $\Gamma$ and we can let $\gamma(a, b)= \{a=b\}$.

  If $h({\cal D})>0$  and ${\cal D}$ ends with a $\hbox{Rep}_1^{r=_r}$-inference, i.e it is of the form:
  
  \[
  \ba{c}
  {\cal D}_0\\
  \underline{a=b , \Gamma^- \seq c=a}\\
  a=b, \Gamma^-\seq c=b
  \ea
  \]
  by induction hypothesis we have that $a=b, \Gamma^-$ is of the form $\gamma'(c, a), \Gamma^{--}$.
  If $a\approx b$ does not belong to  $\gamma'(c, a)$ it suffices to let
   $\gamma(a,b)= \gamma'(c,a)\cup\{ a=b \}$. Otherwise, since  $\gamma'(c, a)$ can be represented as
   \[
   c\approx a_1, \ldots, a_i\approx a, a\approx b, b\approx a_{i+3}, \ldots, a_{n-1}\approx a
   \]
  we can let $\gamma(c,b)= \{c\approx a_1, \ldots, a_i\approx a, a\approx b\}$.
  The same argument applies if  ${\cal D}$ ends with a $\hbox{Rep}_2^{r=_r}$-inference.
  
  b)  If $h({\cal D})=0$ then  $A$ occurs in $\Gamma$  and the claim  holds with $m=0$.
  
  If $h({\cal D})>0$  and ${\cal D}$ ends with a $\hbox{Rep}_1^{r}$-inference, assuming, for the sake of notational  simplicity, that the induction hypothesis holds with $m'=2$, the last inference of ${\cal D}$ has one of the following three forms:
  \[
  \ba{lc}
  i)&~~~\underline{b_1 = b, ~A^\circ[x_1/a_1, x_2/a_2], \gamma'_1(a_1,b_1), \gamma'_2(a_2, b_2), \Gamma^-\seq A^\circ[x_1/b_1, x_2/b_2]}\\
&~~b_1 = b, ~A^\circ[x_1/a_1, x_2/a_2], \gamma'_1(a_1,b_1), \gamma'_2(a_2, b_2), \Gamma^-\seq A^\circ[x_1/b, x _2/b_2]
   \ea
  \]

   \[
  \ba{lc}
  ii)&~~~\underline{b_2 = b, ~A^\circ[x_1/a_1, x_2/a_2], \gamma'_1(a_1,b_1), \gamma'_2(a_2, b_2), \Gamma^-\seq A^\circ[x_1/b_1, x_2/b_2]}\\
 & ~~b_2 = b, ~A^\circ[x_1/a_1, x_2/a_2], \gamma'_1(a_1,b_1), \gamma'_2(a_2, b_2), \Gamma^-\seq A^\circ[x_1/b_1, x _2/b]
   \ea
  \]

   \[
  \ba{lc}
  iii)&~~~\underline{a = b, ~A^\circ[x_1/a_1, x_2/a_2, x/a], \gamma'_1(a_1,b_1), \gamma'_2(a_2, b_2), \Gamma^-\seq A^\circ[x_1/b_1, x_2/b_2, x/a]}\\
  &~~a = b, ~A^\circ[x_1/a_1, x_2/a_2, x/a ], \gamma'_1(a_1,b_1), \gamma'_2(a_2, b_2), \Gamma^-\seq A^\circ[x_1/b_1, x _2/b_2, x/b]
   \ea
  \]

  In case $i)$, if $b_1\approx b$ does not belong to $\gamma'_1(a_1,b_1)$ it suffices to let $\gamma_1(a_1, b)= \gamma'(a_1, b_1)\cup \{ b_1=b\}$
  while if $b_1\approx b$ does belong  to $\gamma'_1(a_1,b_1)$, as in the similar case concerning  $a)$,  we have that  $\gamma'_1(a_1,b_1)$ already contains
  a chain connecting $a_1$ and $b$ that can be taken as $\gamma_1(a_1, b)$. In both cases we let $\gamma_2=\gamma'_2$ so that  $m=m'$.
  
  Case $ii)$ is entirely similar to Case $i)$.  
  
  Finally in Case $iii)$ it suffices to let $\gamma_1=\gamma'_1$, $~\gamma_2=\gamma'_2 $ and $\gamma_3= \{a=b\}$ so that $m=3$. $\Box$

  \

  As an immediate consequence of the two previous lemmas and the admissibility of left weakening  we have the following:

  \begin{proposition}
For languages without function symbols,  a sequent derivable in ${\cal R}_{12}^r$ is derivable also in ${\cal R}_2^{rl}$.

\end{proposition}

 \begin{theorem} \label{no function symbol}
 For languages without function symbols, ${\cal R}_2^{rl}$ is equivalent to ${\cal R}_{12}^r$, hence 
 the structural rules are admissible in ${\bf G3[mic]}^{{\cal R}_2^{rl}}$.
 \end{theorem}
 
 {\bf Proof} By Corollary \ref{sistema completo}  ${\cal R}_2^{rl}$ is a subsystem of ${\cal R}_{12}^r$ and the converse holds by the previous Proposition.
  $\Box$

 \


In the classical case, interpreted in terms of the tableau system introduced in  \cite{J67}, which deals with languages without  function symbols, 
this results is a remarkable improvement of the result in 5.1 of \cite{J67}, 
since it states that not only strictnesss can be required  but  also  that replacement can be restricted to left-right replacement.

\section{Orienting replacement in languages with function symbols}

Let   $\hbox{Rep}_1^{l+}$ and $\hbox{Rep}_2^{l+}$  be the rules  $\hbox{Rep}_1^{l}$ and $\hbox{Rep}_2^{l}$  whose 
instances concerning
equalities  ($E$) are replaced by:

\[
\ba{clccl}
\underline{s=r, E[x/s],E[x/r], \Gamma \seq \Delta}&&~~~~\mbox{and}~~~~&\underline{r=s,~ E[x/s], E[x/r], \Gamma \seq \Delta}&\\
s=r, E[x/s], \Gamma \seq \Delta&&~~&r=s,~ E[x/s], \Gamma \seq \Delta&
\ea
\] respectively.
 
 \ 
 
 Note that, thanks to the admissibility of left weakening,  $\hbox{Rep}_1^{l+}$ and $\hbox{Rep}_2^{l+}$ are strengthening  of $\hbox{Rep}_1^{l}$ and $\hbox{Rep}_2^{l}$ respectively. On the other hand, it is straightforward that  Proposition \ref{Equivalence} extends to such rules as well.

 \begin{proposition}\label{Admissibility of Rep1r}
 The rule ${\em Rep}_1^{r}$ is admissible in ${\cal R}_2^{rl+} = \{ \overline{{\em Ref}}, {\em Rep}_2^{l+}, ~ {\em Rep}_2^{r }\}$.
The same holds with $1$ and $2$ exchanged.
 \end{proposition}
  {\bf Proof}   We may assume that all the rules under consideration replace exactly one occurrence of a term by another (see \cite{PP21} and \cite{PP20a}).
  Then we proceed by induction on the height  $h({\cal D})$ of a  derivation ${\cal D}$ in  $\{ \overline{\hbox{Ref}},\hbox{Rep}_1^r,  \hbox{Rep}_2^{l+}, ~ \hbox{Rep}_2^r \}$ that ends with an $\hbox{Rep}_1^r$-inference and contains no other $\hbox{Rep}_1^r$-inference, to show that ${\cal D}$ can be transformed into a derivation ${\cal D}'$  in  ${\cal R}_2^{rl+}$
    of the same endsequent.
  If $h({\cal D})=1$, then ${\cal D}$ has the form:

  \[
  \ba{c}
 \underline{ r=s, \Gamma\seq \Delta, P[x/r]}\\
  r=s, \Gamma\seq \Delta, P[x/s]
  \ea
  \]
  where $ r=s, \Gamma\seq \Delta, P[x/r]$ is  either  an initial sequent or an instance of $  \overline{\hbox{Ref}}$.
  Case 1. $ r=s, \Gamma\seq \Delta, P[x/r]$ is    an initial sequent. Then we have the following subcases:
  
  Case 1.1.  $(r=s, \Gamma)\cap \Delta\neq \es$,  then  $r=s, \Gamma\seq \Delta, P[x/s]$ is also an initial sequent.
  
  Case 1.2.    $ r=s, \Gamma\seq \Delta, P[x/r]$ is  of the form $r=s,~P[x/r],~\Gamma' \seq \Delta, P[x/r]$. Then 
   ${\cal D}$ can be transformed into:
  
 \[
  \ba{cl}
 \underline{ r=s,~ P[x/s], \Gamma'\seq \Delta, P[x/s]} & \vbox to 0pt{$\hbox{Rep}_2^l$}\\
  r=s,~ P[x/r], \Gamma'\seq \Delta, P[x/s]
  \ea
  \]
Case 1.3. $ r=s, \Gamma\seq \Delta, P[x/r]$ is of the form $r=s, \Gamma\seq \Delta, r=s$. 

  Case 1.3.1.$P\equiv x=s$, hence ${\cal D}$ has  the form:
  
  \[
  \ba{c}
  \underline{r=s, ~\Gamma\seq \Delta, (x=s)[x/r]}\\
  r=s, ~~\Gamma\seq \Delta, (x=s)[x/s]
  \ea
  \]then the  conclusion of ${\cal D}$ is an instance of $\overline{\hbox{Ref}}$, that can be taken as ${\cal D}'$.
  
  Case 1.3.2. $P\equiv s^\circ$, with $s^\circ[x/r]\equiv s$, hence  ${\cal D}$ has the form:
  
  \[
  \ba{c}
r=s^\circ[x/r], \Gamma \seq \Delta,~ r=s^\circ[x/r] \\
\overline{r=s^\circ[x/r], \Gamma \seq \Delta, ~r=s^\circ[x/s^\circ[x/r]]}
    \ea
  \]Then ${\cal D}$ can be transformed into:
  
  \[
  \ba{cl}
 \underline{ r=s^\circ[x/r], ~\Gamma\seq \Delta,~ s^\circ[x/s^\circ[x/r]] = s^\circ[x/s^\circ[x/r]]} &  \vbox to 0pt{$\hbox{Rep}_2^{r}$}\\
  \underline{ r=s^\circ[x/r], ~\Gamma\seq \Delta, ~s^\circ[x/r] = s^\circ[x/s^\circ[x/r]]}&  \vbox to 0pt{$\hbox{Rep}_2^{r}$}\\
r=s^\circ[x/r], ~\Gamma\seq \Delta, ~r  = s^\circ[x/s^\circ[x/r]]
   \ea
  \]

Case 2.  $ r=s, \Gamma\seq \Delta, P[x/r]$ is an instance of $ \overline{\hbox{Ref}}$. Then we have the following subcases:

Case 2.1. The  principal formula is  in $\Delta$. Then 
  $ r=s, \Gamma\seq \Delta, P[x/s]$ is also an instance of $\overline{\hbox{Ref}}$.
  
  Case 2.2. The principal formula is $P[x/r]$. Then
  $P[x/r]$ has the form $t=t$, hence $P$ has the form $t^\circ =t$, or $t=t^\circ$ with $t\equiv t^\circ[x/r]$

  Case 2.2.1. $P\equiv t^\circ =t$. Then  ${\cal D}$ is transformed into:
  
  \[
  \ba{cl}
  \underline{r=s, \Gamma\seq \Delta, t^\circ[x/s]=t^\circ[x/s]}& \vbox to 0pt{$\hbox{Rep}_2^r$}\\
  r=s, \Gamma\seq \Delta, t^\circ[x/s]=t &
  \ea
  \]
  Case 2.2.2.  $P\equiv t=t^\circ$.  Then  ${\cal D}$ is transformed into:
   \[
  \ba{cl}
  \underline{r=s, \Gamma\seq \Delta, t^\circ[x/s]=t^\circ[x/s]}& \vbox to 0pt{$\hbox{Rep}_2^r$}\\
  r=s, \Gamma\seq \Delta, t = t^\circ[x/s]&
  \ea
  \]
  
  If $h({\cal D})>0$ we distinguish the following cases:
  
  Case 3.  The last inference  of the immediate subderivation of ${\cal D}$ is an $\hbox{Rep}_2^r $-inference.
  
  Case 3.1. 
  \[
  \ba{cl}
 \underline{ q=p, ~r=s,~\Gamma\seq \Delta', ~Q[y/p],~P[x/r]}&\vbox to 0pt{$\hbox{Rep}_2^r$}\\
 \underline{ q=p, ~r=s,~\Gamma\seq \Delta', ~Q[y/q],~P[x/r]}&\\
 q=p, ~r=s,~\Gamma\seq \Delta', ~Q[y/q],~P[x/s]
  \ea
  \]
  is transformed into:
  
   \[
  \ba{cl}
 \underline{ q=p, ~r=s,~\Gamma\seq \Delta', ~Q[y/p],~P[x/r]}&\vbox to 0pt{$\hbox{ind}$}\\
 \underline{ q=p, ~r=s,~\Gamma\seq \Delta', ~Q[y/p],~P[x/s]}& \vbox to 0pt{$\hbox{Rep}_2^r$}\\
 q=p, ~r=s,~\Gamma\seq \Delta', ~Q[y/q],~P[x/s]
  \ea
  \]
      
  Case 3.2.
  \[
  \ba{cl}
 \underline{ q=p, ~r=s,~\Gamma\seq \Delta,~P[y/p, x/r]}&\vbox to 0pt{$\hbox{Rep}_2^r$}\\
 \underline{ q=p, ~r=s,~\Gamma\seq \Delta~P[ y/q, x/r]}&\\
 q=p, ~r=s,~\Gamma\seq \Delta,~P[y/q, x/s]
  \ea
  \]
  is transformed into:
  
   \[
  \ba{cl}
 \underline{ q=p, ~r=s,~\Gamma\seq \Delta,~P[y/p, x/r]}& \vbox to 0pt{$\hbox{ind}$}\\
 \underline{ q=p, ~r=s,~\Gamma\seq \Delta~P[ y/p, x/s]}&\vbox to 0pt{$\hbox{Rep}_2^r$}\\
 q=p, ~r=s,~\Gamma\seq \Delta,~P[y/q, x/s]
  \ea
  \]

  Case 3.3.
  \[
  \ba{cl}
  \underline{q^\circ[y/r]=p,~r=s,~\Gamma\seq \Delta,~ P[x/p]}&\vbox to 0pt{$\hbox{Rep}_2^r$}\\
  \underline{q^\circ[y/r]=p,~r=s,~\Gamma\seq \Delta,~ P[x/ q^\circ[y/r]]}&\\
  q^\circ[y/r]=p,~r=s,~\Gamma\seq \Delta,~ P[x/ q^\circ[y/s]]
  \ea
  \]
  is transformed into:
  \[
  \ba{cl}
  \underline{q^\circ[y/r]=p,~r=s,~\Gamma\seq \Delta,~ P[x/p]}&\vbox to 0pt{$\hbox{LW}$}\\
  \underline{q^\circ[y/r]=p,~q^\circ[y/s]=p,~r=s,~\Gamma\seq \Delta,~ P[x/ p]}&  \vbox to 0pt{$\hbox{Rep}_2^r$}\\\
  \underline{q^\circ[y/r]=p, ~q^\circ[y/s]=p, ~r=s,~\Gamma\seq \Delta,~ P[x/ q^\circ[y/s]]}&\vbox to 0pt{$\hbox{Rep}_2^{l+}$}\\
   q^\circ[y/r]=p,~r=s,~\Gamma\seq \Delta,~ P[x/ q^\circ[y/s]]
  \ea
  \]
  
  Case 3.4.
  
  \[
  \ba{cl}
  \underline{q=p,~ r^\circ[x/q]=s, ~\Gamma \seq \Delta, P[x/r^\circ[y/p]]}&\vbox to 0pt{$\hbox{Rep}_2^r$}\\
   \underline{q=p, ~r^\circ[x/q]=s, ~\Gamma \seq \Delta, P[x/r^\circ[y/q]]}&\\
   q=p,~ r^\circ[x/q]=s,~ \Gamma \seq \Delta, P[x/s]
   \ea
  \]
  is transformed into:
  
  \[ 
  \ba{cl}
  \underline{q=p,~ r^\circ[x/q]=s, ~\Gamma \seq \Delta, P[x/r^\circ[y/p]]}& \vbox to 0pt{$\hbox{LW}$}\\
    \underline{q=p, ~r^\circ[x/q]=s,~r^\circ[x/p]=s,~ \Gamma \seq \Delta, P[x/r^\circ[y/p]]}& \vbox to 0pt{$\hbox{ind}$}\\
     \underline{q=p,~ r^\circ[x/q]=s,~r^\circ[x/p]=s,~ \Gamma \seq \Delta, P[x/s]}&\vbox to 0pt{$\hbox{Rep}_2^{l+}$}\\
   q=p, ~r^\circ[x/q]=s, ~\Gamma \seq \Delta, P[x/s]
   \ea
  \]
  Case 4.   The last inference  of the immediate subderivation of ${\cal D}$ is an $\hbox{Rep}_2^{l+}$-inference acting on a formula $Q$ that is not an equality, namely  an $\hbox{Rep}_2^{l}$-inference.

  \[
  \ba{cl}
  \underline{q=p,~ r=s, Q[y/p], \Gamma' \seq \Delta, P[x/r]} &\vbox to 0pt{$\hbox{Rep}_2^l$}\\
  \underline{q=p,~ r=s, Q[y/q], \Gamma' \seq \Delta, P[x/r]}&\\
  q=p, ~r=s, ~Q[y/q], \Gamma' \seq \Delta, P[x/s]
  \ea
  \]
  is transformed into:
  
  \[
  \ba{cl}
  \underline{q=p, r=s, Q[y/p], \Gamma' \seq \Delta, P[x/r]}  &  \vbox to 0pt{$\hbox{ind}$}\\\
  \underline{q=p, r=s, Q[y/p], \Gamma' \seq \Delta, P[x/s]}& \vbox to 0pt{$\hbox{Rep}_2^l$}\\
  q=p, r=s, Q[y/q], \Gamma' \seq \Delta, P[x/s]
  \ea
  \]

  Case 5.  The last inference of the immediate subderivation of ${\cal D}$ is a $\hbox{Rep}_2^{l+}$-inference acting on an equality $E$. 
   In this case we can proceed as in Case 4, by first inverting the last $\hbox{Rep}_1^r$- inference
  with the preceding   $\hbox{Rep}_2^{l+}$-inference and then applying the induction hypothesis.
  $\Box$

 \begin{theorem}\label{Admissibility2}
    The systems  ${\cal R}_{12}^r$ and $ {\cal R}_2^{rl+}$ are equivalent, hence 
  the structural rules are admissible in ${\bf G3[mic]}^{{\cal R}_2^{rl+}}$. The same holds for  ${\cal R}_1^{rl+}=\{ \overline{{\em Ref}},   {\em Rep}_1^{l+}, ~ {\em Rep}_1^{r} \}$. 
    \end{theorem}
  
  {\bf Proof}  Since, by the previous Proposition,  $\hbox{Rep}_1^r$ is admissible in ${\cal R}_2^{rl+}$, ${\cal R}_{12}^r$ is a subsystem of
  ${\cal R}_2^{rl+}$. By Theorem \ref{AdmissibilityR} and Proposition \ref{Equivalence} we have the converse inclusion. $\Box $

 \

Let  ${\cal R}_1^{rl}$ and  ${\cal R}_2^{rl}$
 be $  \{ \overline{\hbox{ Ref}}, \hbox{Rep}_1^{l}, ~ \hbox{Rep}_1^{r }\}$  and $ \{ \overline{\hbox{Ref} }, \hbox{Rep}_2^{l}, ~ \hbox{ Rep}_2^{r }\}$  respectively.
 
 \begin{proposition}\label{RLC}
 ${\cal R}_1^{rl+}$  and ${\cal R}_2^{rl+}$ are  equivalent to ${\cal R}_1^{rl} + \em{LC}^=$  and  ${\cal R}_2^{rl} + \em{LC}^=$  respectively.
 \end{proposition}
 
 {\bf Proof}
$\hbox{Rep}_1^{l+}$ and $\hbox{ Rep}_2^{l+}$ are immediately derivable by means of $\hbox{LC}^=$ from $\hbox{Rep}_1^{l}$ and $\hbox{ Rep}_2^{l}$ respectively. On the other hand 
$\hbox{LC}^=$ is admissible in both  ${\cal R}_1^{rl+}$  and ${\cal R}_2^{rl+}$ by the previous Theorem. $\Box$
 
 \

This naturally leads to what we consider a quite significant  problems left open by our investigation:

 \

 {\bf Question}  Is it possible to extend Theorem \ref{no function symbol}
  to languages endowed with function symbols, namely to replace    ${\cal R}_2^{rl+}$   by   $ {\cal R}_2^{rl}$  in Theorem \ref{Admissibility2}?
  
  \

In the classical case, a positive answer, interpreted  in terms of the alternate tableau system in \cite{F96},
 would mean that  it is possible to require both strictness and restrict replacement to left-right replacement 
 provided the latter is allowed on all atomic and negation of atomic formulae.

\section{Admissibility of the structural rules in systems based on the Left Reflexivity Rule }

As noticed in \cite{PP20},  it is easy to check that  all the structural  rules are admissible in $\{\hbox{Ref}, \hbox{Rep} \}$, so that by Theorem \ref{basic}
we have the following admissibility result, that can be established also by the method in \cite{NvP98} (see Sec. 4, in \cite{PP20a} for full details):

\begin{theorem}\label{theorem1} 
The structural rules are admissible in ${\bf G3[mic]}^{\{ \em{Ref, Rep} \}}$
\end{theorem}

This result can be improved as follows:

 \begin{theorem}\label{theorem2}
The structural rules are admissible in  
${\bf G3[mic]}^{ \{ \em{Ref}, Rep_2^l\}}$.
Therefore ${\bf G3[mic]}^{\{ \em{Ref, Rep} \}}$ and ${\bf G3[mic]}^{ \{ \em{Ref}, Rep_2^l\}}$ are equivalent.
 The same holds for   ${\bf G3[mic]}^{ \{ \em{Ref}, Rep_1^l\}}$.
\end{theorem}
 
 {\bf Proof} As shown in  \cite{PP20a}, $\hbox{Rep}$ is admissible in $\{ \hbox{Ref}, \hbox{Rep}_2^l \}$,  and  $\hbox{Rep}_2^l $ is derivable from $\hbox{Rep} $ by LW.
Therefore $\{\hbox{Ref}, \hbox{Rep}_2^l \}$  and 
 $\{\hbox{Ref}, \hbox{Rep} \}$  are equivalent so that  the first part  follows by Theorem \ref{theorem1}.  
 
 For the  second part we note that, because of the derivability results in  the proof of Proposition \ref{Equivalence}, the Left Symmetry Rule  is admissible in the four systems considered. $\Box$

 \

As it is proved in \cite{F96} the tableau system corresponding to the rules Ref and Rep, namely to the system ${\bf G3[mic]}^{\{ \hbox{Ref}, \hbox{Rep} \}}$ is complete.
Therefore all the tableau systems corresponding to sequent calculi equivalent to ${\bf G3[mic]}^{\{ \hbox{Ref}, \hbox{Rep} \}}$ are complete. 
In particular by Theorem \ref{theorem2} that applies to ${\bf G3[mic]}^{ \{ \hbox{Ref}, \hbox{Rep}_2^l\}}$, which means   that in the tableau system in \cite{F96} pg. 289  strictness can be required without loosing the completeness of the system.

   \section{Counterexamples to the admissibility of the structural rules}
  Since  the weakening rules and the right contraction rule are   admissible in all  the systems consisting of $\overline{\hbox{Ref}}$ and some of the equality rules,  we will 
   concentrate on the possible failure of  the left contraction LC  and/or the  Cut rule.
    By Proposition \ref{Equivalence} and Theorem \ref{AdmissibilityR}, all the axioms and rules for equality that 
  we have considered are admissible in ${\cal R}_{12}^r$. Thus, by Corollary \ref{equivalentsystems}, to show that at least one among LC  and  Cut
  is  not  admissible in  a system ${\cal S}$ 
 it suffices to find a sequent  derivable in ${\cal R}_{12}^r$ 
  but not in ${\cal S}$.  A case of this kind in which  LC is present, thus obviously admissible, and, therefore, Cut is not admissible, is provided by 
   ${\cal S}_1 =\{\overline{\hbox{Ref}},  \hbox{LC},  \hbox{Rep}_2^{l+}, {\hbox{Rep}_1^r}\}$. In fact 
    for $a, b$ and $c$ distinct, the sequent  $ a=c, b=c  \seq a=b$, which is derivable in ${\cal  R}_{12}^r$,
          is  not derivable in ${\cal S}_1$.  As a matter of fact no sequent of the form
            \[
    *)~~~ a=c,\ldots a=c,  b=c, \ldots, b=c, c=c, \ldots, c=c \seq a=b
     \]is derivable in  ${\cal S}_1$, since it can be the conclusion of  LC,  $\hbox{Rep}_2^{l+}$ or $ {\hbox{Rep}_1^r}$-inference
     only if its premiss has already the form $*)$  and no  initial sequent or  instance of $\overline{\hbox{Ref}}$ has that form.
     Clearly the same holds if in ${\cal S}_1$, $\hbox{Rep}_2^{l+}$  is replaced by the more extended rule $\hbox{Rep}$.
      A similar argument applies to      
      ${\cal S}_2 = \{\overline{\hbox{Ref}}, \hbox{LC},  {\hbox{Rep}_1^{l+}}, {\hbox{Rep}_2^r}\}$     
    with respect to the sequent
      $c=b,~ c=a  \seq a=b$ which is derivable in ${\cal R}_{12}^r$ but not in ${\cal S}_2$ and to the system obtained by replacing  ${\hbox{Rep}_1^{l+}}$
      by $\hbox{Rep}'$.   
           While  for the above systems it is the admissibility of  Cut  that fails, 
                       $\{\overline{\hbox{Ref}},  \hbox{Cut}, =_1, =_2\}$  
 is a system in which it is the
    admissibility   of $\hbox{LC}$, actually of   $\hbox{LC}^= $, that fails,   since,   $a=f(a), a=f(a)\seq a=f(f(a))$
     is  derivable, but   $a=f(a)\seq a=f(f(a))$ is not.
 Another example of the same sort  is provided by  $\{\overline{\hbox{Ref}},  \hbox{Cut}, \hbox{CNG}\}$, which is easily seen to be  equivalent to  $\{\overline{\hbox{Ref}},  \hbox{Cut}, =_1, =_2\}$.  
   Although  in general it may happen for a rule not to be admissible in a system but  admissible in a weaker  one, 
     for the system  we are considering,  since the failure of the admissibility of some of the structural rules is witnessed by the underivability of some sequent, which is obviously preserved by weakening a system,
      if they  are not all admissible
      in ${\cal S}$ and ${\cal S}'$ is a subsystem of ${\cal S}$, 
           then they are not all admissible in ${\cal S'}$ either.
        For example, since  $\{\overline{\hbox{Ref}}, \hbox{CNG}\}$ is a subsystem $\{\overline{\hbox{Ref}},  \hbox{Cut}, \hbox{CNG}\}$,  LC and Cut are not both   admissible also  in $\{\overline{\hbox{Ref}}, \hbox{CNG}\}$.
   Actually that  is still  a case in which it is LC to be not admissible, since  Cut remains admissible as it 
   can be easily verified proceeding by induction on the height of the derivation  in $\{\overline{\hbox{Ref}}, \hbox{CNG}\}$  of its  second premiss.
   But
    note that,  by  4) in Proposition \ref{Equivalence} and the analogue for $\hbox{Rep}_2^r$ in the proof of  Theorem \ref{AdmissibilityR},
    it suffices to add to  $\{ \overline{\hbox{Ref}},  \hbox{CNG}\}$ the left contraction rule restricted to equalities $\hbox{LC}^=$
to obtain a system equivalent to  ${\cal R}_{12}^r$
and, therefore,  the admissibility of both LC and Cut.

 \section{Semishortening derivations}
 
 Let us recall from \cite{PP21} the following definition:

 \begin{definition}  Let $\prec$ be any   antisymmetric relation on terms.
An application of an equality rule with operating equality $r=s$ or $s=r$   is said to be {\em nonlengthening} if $s\not\prec r$ and {\em shorthening}  if $r\prec s$. A derivation is said to be {\em semishortening}  if all its equality inferences with index $2$ are nonlengthening and those  with index $1$ are shortening.
 \end{definition}
 The results in \cite{PP21} can be easily adapted to the following context yielding the following Proposition and Theorem:
   
  \begin{proposition} \label{Semishortening} 
 
     If $\Gamma\seq \Delta$ is derivable in ${\cal R}_{12}^r$, then  $\Gamma\seq \Delta$  has a semishortening derivation in ${\cal R}_{12}^{rl+}$.
  \end{proposition}
 
   {\bf Proof} It suffices to show that $\hbox{Rep}_1^{r}$ and $\hbox{Rep}_2^{r} $ are admissible in the calculus  ${\cal R}_{12\prec}^{rl+}$,  namely ${\cal R}_{12}^{rl+}$ with  the applications of  $  \hbox{Rep}_1^{l+}$ and  
  $ \hbox{Rep}_1^{r}  $   required to be shortening, denoted by $  \hbox{Rep}_{1\prec}^{l+}$ and  
  $ \hbox{Rep}_{1\prec}^{r}  $,   and the applications of  $
   \hbox{Rep}_2^{l +}$ and $ \hbox{Rep}_2^{r} $  to be nonlengthening, denoted by  $\hbox{Rep}_{2\prec}^{l +}$ and $ \hbox{Rep}_{2\prec}^{r} $.
  
   We proceed by   induction on the height of a derivation in ${\cal R}_{12\prec}^{rl+}$ of the premiss of a non shortening  $\hbox{Rep}_1^{r}$-inference or  of a 
   lenghtening     $\hbox{Rep}_2^{r}$-inference.
   
   As for a  non shortening  $\hbox{Rep}_1^{r}$-inference, if the derivation of the premiss is an initial sequent or an instance of $\overline{\hbox{Ref}}$
   or ends with a  $\hbox{Rep}_{2\prec}^{r}$ or a  $\hbox{Rep}_{2\prec}^{l+}$ we apply the same transformations used in the proof of Proposition \ref{Admissibility of Rep1r}. Inspection of the various cases  reveals  that in  the transformed derivation, the given non shortening  $\hbox{Rep}_1^{r}$-inference is replaced by a  $\hbox{Rep}_2^{l}$-inference that,  having  the same operating equality,  turns out to be non lenghtening. Furthermore if the derivation of the premiss ends with a 
    $\hbox{Rep}_{1\prec}^{r}$ or a  $\hbox{Rep}_{1\prec}^{l+}$-inference we can perform similar tranformations leading to a derivation in ${\cal R}_{12\prec}^{rl+}$
    of the conclusion.
    The case of a lenghtening  $\hbox{Rep}_2^{r}$-inference is dealt with in a similar way. We leave the  details to the reader.
    $\Box$

   \begin{theorem}
  The systems ${\cal R}_{12}^r$ and ${\cal R}_{12\prec}^{rl+}$ are equivalent, hence the 
   structural rules are admissible in ${\bf G3[mic]}^{{\cal R}_{12\prec}^{rl+}}$.
  
   \end{theorem}
   
  {\bf Proof} By the previous Proposition, ${\cal R}_{12}^r$ is a subsystem of
  ${\cal R}_2^{rl+}$. The conclusion follows by Theorem \ref{AdmissibilityR} and Proposition \ref{Equivalence} $\Box$

  \

The proof  of Proposition \ref{Semishortening}  uses  the strengthened form  $ \hbox{Rep}_1^{l+}, 
   \hbox{Rep}_2^{l+}$ of  the rules $\hbox{Rep}_1^{l}, 
   \hbox{Rep}_2^{l}$.
      However we have no counterexample, i.e. no particular $\prec$, showing that Proposition \ref{Semishortening} does not hold    for ${\cal R}_{12\prec }^{rl}$, in particular, according to the problem at the end of Section 7, since  ${\cal R}_{12\es}^{rl+}$ amounts to  the same as ${\cal R}_{2}^{rl+}$, we do not have one for  $\prec =\es$.

      \ 
      
      {\bf Note} 
   In case  $\prec$ is the relation induced by $rank$-comparison i.e.  if $r\prec s$ if and only if
   the height (of the  formation tree) of $r$ is smaller than that of $s$, the derivability in ${\cal R}_{12\prec}^{rl}$ is  closely related to the notion of a sequent being 
    {\em directly demonstrable} as defined and claimed to be decidable in  \cite{K63}, pg.90.

 \section{Conclusion}
 
 We have shown how the Gentzen's sequent calculi for first order logic with equality studied in \cite{PP21} naturally evolve into their structural free counterparts
 based on Dragalin's multisuccedent calculi for minimal, intuitionistic  and classical logic.
 From the historical point of view it is worth mentioning that, in the classical case, the system based on ${\cal R}_{12}^r$, that we regard as the most natural one, is the system introduced and semantically investigated in the classic \cite{W60}.  We have shown that various restrictions limiting the scope of the replacement in the equality rules leave all the structural rules admissible. In the classical case all such results ensure  the possibility of placing corresponding restrictions on the semantic tableau method for
 first order logic with equality. A particularly significant result is the possibility of  imposing  strictness as well as orientation of the replacement of equals  
 in case the language lacks function symbols.
On the way of extending this orientability result to general languages we have shown its  reducibility to   the admissibility of the Left Contraction Rule for equalities.
 Whether or not orientability can be obtained without adding such a contraction rule remains an open problem to be settled. 
 Furthermore  we have discussed to what extent   the results in \cite{PP21} concerning semishortening derivations can be extended to the present context leaving open a problem that includes
 the previous one as a particular case.

  \end{document}